\documentclass[11pt]{amsart}
\usepackage{amssymb}
\usepackage{amsthm}
\usepackage{enumerate}

\renewcommand{\phi}{\varphi}

\newcommand{\Kth}{K_\theta}
\newcommand{\kp}{K^p_\theta}
\newcommand{\T}{\mathbb T}
\newcommand{\R}{\mathbb R}

\newcommand{\Z}{\mathbb Z}
\newcommand{\D}{\mathbb D}

\newcommand{\tz}{\mathbb T}
\newcommand{\dd}{\mathbb D}
\newcommand{\vep}{\varepsilon}

\newcommand{\C}{\mathcal C}
\newcommand{\Dd}{\mathcal D}
\newcommand{\Tau}{\mathcal T}

\renewcommand{\le}{\leqslant}
\renewcommand{\ge}{\geqslant}
\newcommand{\x}{\!\cdot\!}
\newcommand{\PW}{{\mathcal PW}}
\newcommand{\Alambda}{T_\lambda}

\newcommand{\into}{\hookrightarrow}

\renewcommand{\Re}{{\rm Re}\,}

\newtheorem{Thm}{Theorem}[section]
\newtheorem{Lem}[Thm]{Lemma}
\newtheorem{Prop}[Thm]{Proposition}
\newtheorem{Cor}[Thm]{Corollary}

\newtheorem{Ex}[Thm]{Example}
\newtheorem{Conj}[Thm]{Conjecture}

\def\beginpf{\medskip\noindent {\bf Proof.} \;}

\title[Symbols of truncated Toeplitz operators]
{Symbols of truncated Toeplitz operators}

\author{Anton Baranov, Roman Bessonov, Vladimir Kapustin}

\address[A. Baranov]{Department of Mathematics and Mechanics, St. Petersburg State
University, 28, Universitetskii pr., St. Petersburg, 198504, Russia}
\email{anton.d.baranov@gmail.com}
\address[R. Bessonov and V. Kapustin]
{St.Petersburg Department of Steklov Mathematical Institute RAS,
27, Fontanka, St.Petersburg, 191023, Russia}

\email{bessonov@pdmi.ras.ru, kapustin@pdmi.ras.ru}

\date{}

\thanks{The authors are supported by the RFBR grant 08-01-00723.
The first author is supported by grant MK-7656.2010.1 and by
Federal Program of Ministry of Education 2010-1.1-111-128-033.}

\begin{document}

\begin{abstract}
We consider three problems connected with coinvariant subspaces
of the backward shift operator in Hardy spaces $H^p$: 

-- properties of truncated Toeplitz operators;

-- Carleson-type embedding theorems for the coinvariant subspaces;

-- factorizations of pseudocontinuable functions from $H^1$.

\noindent These problems turn out to be closely connected and even, in a sense, 
equivalent. The new approach based on the factorizations 
allows us to answer a number of challenging questions about 
truncated Toeplitz operators posed by Donald Sarason.
\end{abstract}

\maketitle
\sloppy
\section{Introduction}\label{introd}

Let $H^p$, $1\le p\le \infty$, denote the Hardy space in the unit disk $\dd$,
and let $H^p_- = \overline{zH^p}$. 
As usual, we identify the functions in $H^p$ in the 
disk and their nontangential boundary values on the unit circle $\T$. 

A function $\theta$ which is  analytic and bounded in $\dd$
is said to be {\it inner} if
$|\theta| =1$ m-a.e. on $\T$ in the sense of nontangential boundary values;
by $m$ we denote the normalized Lebesgue measure on $\tz$.
With each inner function $\theta$ we associate the subspace
$$
K^p_{\theta}=H^p\cap \theta H^p_-
$$
of the space $H^p$. Equivalently, 
one can define $K_\theta^p$ as the set of all functions in $H^p$
such that $\langle f,\theta g \rangle=
\int_\tz f\, \overline{\theta g} \,  dm=0$
for any $g\in H^q$, $1/p+1/q =1$. In particular, 
$$
K_\theta :=K_\theta^2 = H^2\ominus \theta H^2
$$
(in what follows we usually omit the exponent 2).
It is well known that, for $1 \le p \le \infty$, any closed subspace 
of $H^p$ invariant with respect to the backward shift 
$(S^*f)(z)=\frac{f(z)-f(0)}{z}$
is of the form $K_\theta^p$ for some inner function $\theta$ 
(see \cite[Chapter II]{ga} or \cite{cimros}). 
Subspaces $K_\theta^p$ are often called
{\it star-invariant subspaces}. These subspaces 
play an outstanding role both in function and operator theory
(see \cite{cimros, nik, nk2}) and, in particular,
in the Sz.-Nagy--Foias model for contractions in a Hilbert space
(therefore they are sometimes referred to  as {\it model subspaces}).

A characteristic property of the elements of the spaces $\kp$ is the existence of
a {\it pseudocontinuation} outside the unit disk: if $f\in \kp$, then there exists 
a function $g$ which is meromorphic and of Nevanlinna class  in $\{z: |z|>1\}$
such that $g= f$ a.e. on $\T$ in the sense of nontangential boundary values.
\medskip

Now we discuss in detail the three main themes 
of the paper as indicated in the abstract.

\subsection{Truncated Toeplitz operators on $\Kth$} 
Recall that the classical Toeplitz operator on $H^2$
with a symbol $\phi \in L^\infty(\T)$
is defined by $T_\phi f = P_+(\phi f)$, $f\in H^2$,
where $P_+$ stands for the orthogonal projection from 
$L^2:= L^2(\T, m)$ onto $H^2$.

Now let $\phi\in L^2$. We define the {\it truncated Toeplitz operator} $A_\phi$
on bounded functions from $\Kth$ by the formula 
$$
A_\phi f = P_{\theta} (\phi f), \qquad f\in \Kth\cap L^\infty(\T). 
$$ 
Here $P_\theta f= P_+f - \theta P_+ (\overline \theta f)$ 
is the orthogonal projection onto $\Kth$.
In contrast to the Toeplitz operators on $H^2$ (which satisfy $\|T_\phi \| = 
\|\phi\|_\infty$), the operator $A_\phi$ may be extended to a bounded operator 
on $\Kth$ even for some unbounded symbols $\phi$. 
The class of all bounded truncated Toeplitz 
operators on $\Kth$ will be denoted by $\Tau(\theta)$.
\medskip

Certain special cases of truncated Toeplitz operators are well known and play
a prominent role in the operator theory. If $\phi(z) = z$, then 
$A_\phi = S_\theta$ is the so-called restricted shift operator, the scalar 
model operator from the Sz.-Nagy--Foias theory. If $\phi \in H^\infty$, 
then $A_\phi=\phi(S_\theta)$. Truncated Toeplitz operators include
all finite Toeplitz matrices (corresponding to the case $\theta(z) = z^n$)
and the Wiener--Hopf convolution operators on an interval which are unitary 
equivalent to the truncated Toeplitz operators on the space 
generated by the singular inner function 
associated with a point mass on the circle, 
\begin{equation}
\label{tha}
\theta_a(z) = \exp\Big(a \frac{z+1}{z-1}\Big)
\end{equation}
(for a more detailed discussion see \cite{bar}). 
Bercovici, Foias and Tannenbaum studied truncated (or {\it skew}) Toeplitz 
operators (mainly, with symbols which are rational functions with pole at zero)
in connection with control theory (see \cite{ber1, ber2}). However, a 
systematic study of truncated Toeplitz operators with $L^2$ symbols was started
recently by Sarason in \cite{sar}. This paper laid the basis 
of the theory and inspired much of the subsequent activity in the field
\cite{bar, gar1, gar2, gar3}.
\bigskip

Unlike standard Toeplitz operators on $H^2$, the symbol of a truncated Toeplitz 
operator is not unique. The set of all symbols of an arbitrary operator 
$A_\phi$ is exactly the set $\phi+\theta H^2+\overline{\theta H^2}$, see 
\cite{sar}. Clearly, any bounded function $\phi\in L^\infty$ determines
the bounded operator $A_\phi$ with the norm 
$\|A_\phi\| \leqslant \|\phi\|_{\infty}$. The first basic question 
on truncated Toeplitz operators posed in \cite{sar} is {\it whether 
every bounded operator $A_\phi$ has a bounded symbol}, i.e., is
a restriction of a bounded Toeplitz operator on $H^2$. Note that if 
a  truncated Toeplitz operator with a symbol $\phi \in H^2$ is bounded,
then, as a consequence of the commutant lifting theorem,
it admits an $H^\infty$ symbol (see \cite[Section 4]{sar}). 
On the other hand, by  the results of Rochberg \cite{roch}  
(proved in the context of 
the Wiener--Hopf operators and the Paley--Wiener spaces) any operator 
in $\Tau(\theta_a)$  has a bounded symbol.
However, in general the answer to this question is negative:
in \cite{bar} inner functions $\theta$ are constructed, for which
there exist operators in $\Tau(\theta)$ (even of rank one) that have 
no bounded symbols. 

Thus, the following question seems to be of interest: 
{\it in which spaces $\Kth$ does any bounded truncated Toeplitz operator 
admit a bounded symbol?} In the present paper
we obtain a description of such inner functions. In particular,
we show that this is true for an interesting class of 
{\it one-component inner functions} introduced by Cohn in \cite{cohn82}: 
these are functions $\theta$ such that the sublevel set 
$$
\{z\in \D:\, |\theta(z)| <\vep\}
$$ 
is connected for some $\vep \in (0,1)$. This statement was conjectured 
in \cite{bar}. A basic example of a one-component inner function is 
the function $\theta_a$ given by (\ref{tha}).
\medskip

\subsection{Embeddings of the spaces $K^{p}_{\theta}$}
Let $\mu$ be a finite positive Borel measure in the closed unit disk
$\overline{\D}$. We are interested in the class of measures such that 
Carleson-type embedding $\kp \into L^p(\mu)$ is bounded. 
Since the functions in $\kp$ are well-defined only 
$m$-almost everywhere on $\T$, one should be careful when dealing 
with the restriction of $\mu$ to $\T$. Recall that, by a theorem 
of Aleksandrov, functions in $\kp$ which are continuous in the closed 
disk $\overline{\D}$ are dense in $\kp$, see \cite{al8} or \cite{cimmat}.
(While this statement is trivial for the Blaschke products, there is no
constructive way to prove the statement in the general case.)
This allows one to define the embedding on the dense set of all continuous
functions from $\Kth$ in a natural way
and then ask if it admits a bounded continuation to the whole space $\kp$.
However, this extension may always be viewed as an embedding operator
due to the following theorem by Aleksandrov.

\medskip

\begin{Thm}\cite[Theorem 2]{al1}
\label{ab}  
Let $\theta$ be an inner function, let $\mu$ be a positive 
Borel measure on $\T$, and let $1\le p<\infty$. 
Assume that for any continuous function $f\in \kp$ 
we have 
\begin{equation}
\label{emb1}
\|f\|_{L^p(\mu)} \le C\|f\|_p.
\end{equation} 
Then all functions from $K_\theta^p$ possess angular boundary values $\mu$-almost 
everywhere and for any $f\in K_\theta^p$, 
$(\ref{emb1})$ holds for its boundary values. 
\end{Thm}


The angular convergence $\mu$-almost everywhere gives us a nice illustration
of how the embedding acts. This approach, essentially based on 
results of Poltoratski's paper \cite{polt}, uses deep analytic techniques.
For our purposes we will need the $L^2$-convergence, which can be 
established much simpler. To make the exposition more self-contained 
we present the corresponding arguments in Section~\ref{ele}.

\medskip

Denote by $\Dd_p(\theta)$ the class of finite complex 
Borel measures $\mu$ on 
$\overline{\D}$ for which the embedding $\kp \subset L^p(|\mu|)$
holds; by $|\mu|$ we denote the total variation of the complex measure 
$\mu$. The class of positive measures in $\Dd_p(\theta)$
is denoted by $\Dd_p^+(\theta)$. The classes $\Dd_p^+(\theta)$ 
contain all Carleson
measures, i.e., measures for which the embedding $H^p\into L^p(\mu)$
is a bounded operator (for some, and hence for all $p >0$).
However, the class $\Dd_p(\theta)$ is usually much wider due to additional
analyticity (pseudocontinuability) of the elements of $\kp$ on the boundary.
The problem of the description of the class $\Dd_p(\theta)$
for general $\theta$ was posed by Cohn in 1982; it is still open.
Many partial results may be found in 
\cite{cohn82, cohn86, vt, al2, nv02, bar1, bar2}.
In particular, the classes $\Dd_p(\theta)$ are described if $\theta$ 
is a one-component inner function; in this case there exists a nice geometric
description analogous to the classical Carleson embedding theorem 
\cite{vt, al4, al2}. 
Moreover, Aleksandrov \cite{al2} has shown that $\theta$ is one-component
if and only if all classes $\Dd_p(\theta)$ for $p>0$ coincide. 

In what follows we denote by $\C_p(\theta)$ the set of finite 
complex Borel 
measures $\mu$ on the unit circle $\T$ such that $|\mu|$ 
is in $\Dd_p(\theta)$; the class of positive measures in $\C_p(\theta)$
will be denoted by $\C_p^+(\theta)$.

If $\mu \in \C_2(\theta)$, we may define the bounded operator $A_\mu$ 
on $\Kth$ by the formula
\begin{equation}\label{amu}
(A_\mu f, g)=\int f\bar g\,d\mu.
\end{equation}
It is shown in \cite{sar} that $A_\mu \in \Tau(\theta)$. 
This follows immediately from the following 
characteristic property of truncated Toeplitz operators.

\begin{Thm} \cite[Theorem 8.1]{sar} 
\label{sarason} 
A bounded operator $A$ on $\Kth$ is a truncated Toeplitz 
operator if and only if the condition 
$f, zf\in K_{\theta}$ yields $(Af,f)=(Azf, zf)$.
\end{Thm}

A complex measure $\mu$ on $\T$ with finite total variation
(but not necessarily from $\C_2(\theta)$)
will be called a {\it quasisymbol} for a truncated Toeplitz operator $A$ if
$(Af, g)=\int f\bar g\,d\mu$ holds for all continuous functions $f, g\in\Kth$. 
The symbol $\phi$ of $A=A_\phi$ can be regarded as
a quasisymbol if we identify it with the measure $\phi m$.
 
The following conjecture was formulated by Sarason in \cite{sar}: 
{\it every bounded truncated Toeplitz operator $A$ coincides with $A_\mu$ 
for some $\mu\in \C_2(\theta)$}. Below we prove this conjecture. Moreover, 
we show that nonnegative bounded truncated Toeplitz 
operators are of the form $A_\mu$ with $\mu\in\C_2^+(\theta)$. 
We also prove that truncated Toeplitz operators with bounded symbols 
correspond to complex measures from the subclass $\C_1(\theta^2)$
of $\C_2(\theta^2)=\C_2(\theta)$, and, finally, that every bounded truncated 
Toeplitz operator has a bounded symbol if and only if 
$\C_1(\theta^2)=\C_2(\theta^2)$.

\medskip

\subsection{Factorizations}
Now we consider a factorization problem for pseudocontinuable functions 
in $H^1$, which proves to have an equivalent reformulation in terms of 
truncated Toeplitz operators. 

It is well known that any function $f \in H^1$ can be represented as 
the product of two functions $g, h\in H^2$ with $\|f\|_1 = \|g\|_2\x \|h\|_2$. 
By the definition of the spaces $\kp$, there is a natural involution on 
$\Kth$:
\begin{equation}\label{invo}
f\mapsto \tilde f=\bar z \theta \bar f\in\Kth, \qquad f\in\Kth.
\end{equation}
Hence, if $f,g \in \Kth$, then $fg\in H^1$ and $\bar z^2 
\theta^2 \bar f\bar g \in H^1$. Thus,
$$
fg \in H^1\cap \bar{z}^2 \theta^2 \overline{H^{1}} = H^1 \cap \bar{z}\theta^2 H^{1}_{-}.
$$
If $\theta(0)=0$, then $\theta^2/z$ is
an inner function and $H^1\cap \bar{z}\theta^2 H^{1}_{-}= K^{1}_{\theta^2/z}$. 

It is not difficult to show that linear combinations of products of pairs 
of functions from $\Kth$ form a dense subset of $H^1\cap \bar{z}\theta^2 H^{1}_{-}$. 
We are interested in a stronger property: 

{\it For which $\theta$ may any function 
$f\in H^1 \cap \bar{z}\theta^2 H^{1}_{-}$ be represented in the form}
\begin{equation}
\label{fact1}
f=\sum_k g_kh_k, \qquad  g_k, h_k \in \Kth, \qquad 
\sum_k \|g_k\|_2\x  \|h_k\|_2 <\infty?
\end{equation}
We still use the term {\it factorization} 
for the representations of the form (\ref{fact1}), by analogy with 
the usual row-column product.

Below we will see that, for functions $f\in H^1 \cap \bar{z}\theta^2 H^1_-$, 
not only a usual factorization $f=g\x h$, $g,h\in \Kth$, but even a weaker 
factorization (\ref{fact1}) may be impossible.  
We prove that this problem is equivalent 
to the problem of existence of bounded symbols for all bounded truncated 
Toeplitz operators on $\Kth$.  

Let us consider two special cases  of the problem.
Take $\theta(z)=z^{n+1}$. The spaces $\Kth$ and $K^{1}_{\theta^2/z}$ 
consist of polynomials of degree at most $n$ and $2n$, respectively,
and then, obviously, $K^{1}_{\theta^2/z} = \Kth \cdot \Kth$. 
However, it is not known if a norm controlled factorization
is possible, i.e., if for a polynomial $p$ of degree at most $2n$ 
there exist polynomials $q, r$ of degree at most $n$ such that $p=q\x r$ 
and $\|q\|_2\x\|r\|_2 \le C \| p\|_1$, where $C$ is an absolute constant
independent on $n$. On the other hand, it is shown in \cite{Volberg} that 
there exists a representation $p=\sum_{k=1}^4 q_kr_k$ with 
$\sum_{k=1}^4 \|q_k\|_2\x\|r_k\|_2\le C\|p\|_1$.
 
For $\theta=\theta_a$ defined by (\ref{tha}), 
the corresponding model subspaces $\kp$
are the natural analogs of the Paley--Wiener spaces $\PW^p_a$ of entire 
functions. The space $\PW^p_a$ consists of all entire functions of exponential
type at most $a$, whose restrictions to $\R$ are in $L^p$. 
It follows from our results (and may be proved directly,
see Section \ref{onecomp}) that every entire function $f\in \PW^1_{2a}$ 
of exponential type at most $2a$ and summable on the real line $\R$ admits
a representation $f = \sum_{k=1}^4 g_k h_k $  
with $f_k, g_k\in\PW^2_a$, 
$\sum_{k=1}^{4}\|g_k\|_2\x\|h_k\|_2\le C\|f\|_1$.

\bigskip
\section{Main results}\label{main}

Our first theorem answers Sarason's question 
about representability of bounded truncated
Toeplitz operators via Carleson measures for $\Kth$.

\begin{Thm} \label{t1}
$1)$ Any nonnegative bounded truncated Toeplitz operator on $\Kth$ 
admits a quasisymbol which is a nonnegative measure from $\C_2^+(\theta)$.

$2)$ Any bounded truncated Toeplitz operator on $\Kth$ admits
a quasisymbol from $\C_2(\theta)$.
\end{Thm}

In general, in assertion $1)$ of the theorem, $\mu$ cannot be chosen 
absolutely continuous, i.e., bounded nonnegative truncated Toeplitz 
operators may have no nonnegative symbols. 
Let $\delta$ be the Dirac measure at a point of $\T$, for which 
the reproducing kernel belongs to $\Kth$. Then the operator $A_\delta$ 
cannot be realized by a nonnegative symbol unless the dimension of $\Kth$ is 1.
Indeed, if $\mu$ is a positive absolutely continuous measure, then 
the embedding $\Kth\into L^2(\mu)$ must have trivial kernel, while 
in this example it is a rank-one operator.

\medskip

The next theorem characterizes operators from $\Tau(\theta)$ that have 
bounded symbols.

\begin{Thm} \label{c1}
A bounded truncated Toeplitz operator $A$ admits a bounded symbol if and 
only if $A=A_\mu$ for some $\mu\in\C_1(\theta^2)$. 
\end{Thm}

\medskip

In the proofs of these results the key role is played by the following 
Banach space $X$ defined by
\begin{equation}
\label{spaceX}
X=\bigg\{\sum_k x_k\bar y_k: \quad x_k, y_k\in\Kth, \qquad
\sum_k \|x_k\|_2\x \|y_k\|_2<\infty \bigg\}.
\end{equation}
The norm in $X$ is defined as the infimum of $\sum \|x_k\|_2\x \|y_k\|_2$
over all representations of the element in the form $\sum x_k\bar y_k$.

\begin{Thm} \label{t2}
$1)$ The space dual to $X$ can be naturally identified with $\Tau(\theta)$.
Namely, continuous linear functionals over $X$ are of the form
\begin{equation} 
\label{overx}
   \Phi_A(f) = \sum_k (Ax_k, y_k), \qquad f = \sum_k x_k \bar{y}_k \in X,
\end{equation}
with $A\in\Tau(\theta)$, and the correspondence between the functionals 
over $X$ and the space $\Tau(\theta)$ is one-to-one and isometric.

$2)$ With respect to the duality $(\ref{overx})$ the space $X$ is dual 
to the class of all compact truncated Toeplitz operators. 
\end{Thm}

The next theorem establishes a connection between the factorization problem,
Carleson-type embeddings, and the existence of a bounded symbol for
every bounded truncated Toeplitz operator on $\Kth$.

\begin{Thm} \label{t3}
The following are equivalent:

$1)$\, any bounded truncated Toeplitz operator on $\Kth$ admits a bounded
symbol;

$2)$\, $\C_1(\theta^2)=\C_2(\theta^2)$;

$3)$\, for any $f\in H^1\cap \bar z\theta^2 H^1_-$ there exist
$x_k, y_k\in\Kth$ with $\sum_k \|x_k\|_2\x \|y_k\|_2<\infty$ 
such that $f=\sum_k x_k y_k$. 

\end{Thm}

In the proof it will be shown that condition 2) can be replaced by the
stronger condition 

\medskip

$2')$\, $\Dd_1(\theta^2)=\Dd_2(\theta^2)$.

\medskip

Condition 3) also admits formally stronger, but equivalent
reformulations. If 3) is fulfilled, then, 
by the Closed Graph theorem, one can always find $x_k, y_k$ such that 
$\sum_k \|x_k\|_2\cdot \|y_k\|_2 \le C\|f\|_1$ for some constant $C$ 
independent from $f$. Thus, 3) 
means that $X = H^1\cap \bar z\theta^2 H^1_-$ 
and the norm of $X$ is equivalent to $L^1$ norm. Moreover, 
it follows from Proposition \ref{classx} that one can require 
that the sum contain at most four summands.

\medskip

If $\theta$ is a one-component inner function, then all classes $\C_p(\theta)$
coincide, see \cite[Theorem 1.4]{al2}. If $\theta$ is one-component, then
$\theta^2$ is, too, hence $\C_1(\theta^2)=\C_2(\theta^2)$. As an immediate 
consequence of Theorem \ref{t3}
we obtain the following result conjectured in \cite{bar}:

\begin{Cor}
\label{onecomp1} If $\theta$ is a one-component inner function,
then the equivalent conditions of Theorem $\ref{t3}$
are fulfilled.
\end{Cor}  

We do not know if the converse is true, that is, whether 
the equality $\C_1(\theta^2)=\C_2(\theta^2)$ implies that $\theta$
is one-component. If this is true, it 
would give us a nice 
geometrical description of inner functions $\theta$ 
satisfying the equivalent conditions of Theorem \ref{t3}:

\begin{Conj} 
The equivalent conditions of Theorem $\ref{t3}$ are fulfilled if and only if 
$\theta$ is one-component.
\end{Conj}

Theorem \ref{t3} also allows to extend considerably the class of 
counterexamples to the existence of a bounded symbol.
Let us recall the definition of the Clark measures $\sigma_\alpha$ \cite{cl}.
For each $\alpha\in\mathbb{T}$ there exists a finite (singular) positive 
measure $\sigma_\alpha$  on $\T$ such that
\begin{equation}
\label{clark}
\Re \frac{\alpha+\theta(z)}{\alpha-\theta(z)}=
\int\limits_\mathbb{T} \frac{1-|z|^2}{|1-\bar \tau z|^2}\,
d \sigma_\alpha(\tau), \qquad z\in\mathbb{D}.
\end{equation}
If $\sigma_\alpha$ is purely atomic,
i.e., if $\sigma_\alpha = \sum_n a_n\, \delta_{t_n}$,
then the system $\{k_{t_n}\}$ is an orthogonal basis in ${K^2_{\theta}}$;
in particular, 
$k_{t_n}\in {K^2_{\theta}}$ and $\|k_{t_n}\|_2^2 = |\theta'(t_n)|/(2\pi)$.
\smallskip

It is shown in \cite[Theorem 8]{al1} that the condition
$\C_1(\theta^2) = \C_2(\theta^2)$ implies 
that all measures $\sigma_\alpha$ are discrete. 

\begin{Cor}
\label{onecomp5} 
If, for some $\alpha\in\T$, the Clark measure $\sigma_\alpha$ 
is not discrete, then the conditions of Theorem~$\ref{t3}$ do not hold 
and, in particular, there exist operators 
from $\Tau(\theta)$ that do not admit a bounded symbol.
\end{Cor}  

\bigskip

\section{Embeddings $\Kth\into L^2(\mu)$: the radial $L^2$-convergence}
\label{ele}
In this section we present a more elementary 
approach to embedding theorems which is different from 
that of Theorem~\ref{ab}. Sometimes it may be more convenient to work 
in the $L^2$-convergence 
setting than with continuous functions from $\Kth$. Here we impose 
an extra assumption $\theta(0)=0$, or, equivalently, $1\in\Kth$,
to which the general case can easily be reduced (via
transform (\ref{crof}) defined below), but we omit the details 
of the reduction. 

We will show that the condition 
$\mu\in\C_2^+(\theta)$ is equivalent to the existence of
an operator $J:\Kth\to L^2(\mu)$ such that 

$(i)$ \; if $f,zf\in\Kth$ then $Jzf=zJf$,

$(ii)$ \, $J1=1$. 

\noindent Moreover, these properties uniquely determine the operator, 
which turns out to coincide with the embedding operator
$\Kth\into L^2(\mu)$ defined by Theorem~\ref{ab}. The proofs are based
on the following result of Poltoratski, see also \cite{vk}. 
For $g\in\Kth$, $g_r$ denotes the function $g_r(z)=g(rz)$.

\begin{Prop} (cf. \cite[Theorem 1.1]{polt})
\label{polt} 
Let $\theta(0)=0$. If a bounded operator  $J: \Kth\to L^2(\mu)$ satisfies 
the properties $(i), (ii)$, then for any $g\in\Kth$ we have 
$\|g_r\|_{L^2(\mu)}\le 2\cdot \|Jg\|_{L^2(\mu)}$ 
and $g_r\to Jg$ in $L^2(\mu)$ as $r\nearrow 1$.
\end{Prop}

\beginpf
Consider the Taylor expansion of $g\in\Kth$, 
$$
g(z)=\sum_{k=0}^\infty a_kz^k,
$$
and introduce the functions $g_n\in\Kth$, 
$$
g_n(z)=\sum_{k=0}^\infty a_{k+n}z^k.
$$
By induction from the relation $J g_n=a_n+zJ g_{n+1}$ we obtain the formula 
$$
Jg=\sum_{k=0}^{n-1}a_kz^k+z^nJg_n.
$$
We have $\|g_n\|_2\le \|g\|_2$ and $\|g_n\|_2\to 0$ as $n\to +\infty$. 
Therefore, $\sum_{k=0}^{n-1}a_kz^k\to Jg$ in $L^2(\mu)$. Since $g_r$ are 
the Abel means of the sequence $\Big(\sum_{k=0}^{n-1}a_kz^k\Big)_{n\ge 1}$, 
we conclude that $g_r\to Jg$ as well. \qed

\medskip

Since for a continuous function $g\in\Kth$, $Jg$ coincides with $g$ 
$\mu$-almost everywhere, $J$ is the same operator as the embedding 
from Theorem~\ref{ab}.

\medskip

By Proposition \ref{polt} 
the function $z^{-1}\theta(z)$ (or $\frac{\theta(z)-\theta(0)}{z}$ 
in the general case, if $\theta(0)\ne0$) has the boundary function defined
by the limit in $L^2(\mu)$ of $\theta_r$, where $\theta_r(z)=\theta(rz)$. 
This allows us to define the boundary values of $\theta$ $\mu$-almost 
everywhere.

\begin{Prop} \label{mod1}
If $\mu\in\C_2(\theta)$, then $|\theta|=1$ $|\mu|$-almost everywhere.
\end{Prop}

This fact is mentioned in \cite{al1} and its proof there seems to use 
the techniques of convergence $\mu$-almost everywhere. A more elementary proof 
is given below for the reader's convenience.


\beginpf
We may assume that $\mu\in\C_2^+(\theta)$.
It is easy to check the relation 
$$
M_zJ-JA_z=(\cdot, \bar z\theta)\theta,
$$
where $M_z$ is the operator of multiplication by $z$ on $L^2(\mu)$, 
$J$ is the embedding $\Kth\into L^2(\mu)$, $A_z$ is the truncated Toeplitz 
operator with symbol $z$, i.e., the model contraction $S_\theta$. 
Indeed, on vectors orthogonal to $\bar z\theta$ both sides equal 0, 
and for $\bar z\theta$ the formula can be verified by a simple 
straightforward calculation. Similarly, 
$M_{\bar z}J-JA_{\bar z}=(\cdot, 1)\bar z$, hence 
$$
J^\ast M_z-A_z J^\ast=(\cdot, \bar z)1.
$$

We obtain
$$
\begin{aligned}
JJ^\ast M_z-M_z JJ^\ast=&J(J^\ast M_z-A_z J^\ast)-(M_zJ-JA_z)J^\ast\\
=&(\cdot, \bar z)J1-(\cdot, J\bar z\theta)\theta
=(\cdot, \bar z)1-(\cdot, \bar z\theta)\theta.
\end{aligned}
$$
Theorem~6.1 of \cite{vk-zns} says that if $K=\sum(\cdot, \bar u_k)v_k$
is a finite rank (or even trace class) operator on $L^2(\mu)$, where $\mu$ 
is a singular measure on $\T$, and if $K=XM_z-M_zX$ for
some bounded linear operator $X$ on $L^2(\mu)$, then $\sum u_kv_k=0$ 
$\mu$-almost everywhere. By this theorem $z-z|\theta|^2=0$, hence
$|\theta|=1$ $\mu$-almost everywhere, as required. \qed

\bigskip
\section{The space $X$}\label{x}

As above, the space $X$ is defined by formula (\ref{spaceX}),
$$
X=\left\{\sum x_k\bar y_k: \quad x_k, y_k\in\Kth, \quad
\sum \|x_k\|_2\x \|y_k\|_2<\infty \right\}.
$$
We also consider the analytic analog $X_a$ of the space $X$, 
\begin{equation}
\label{thetaX}
X_a=\left\{\sum x_k y_k: \quad x_k, y_k\in\Kth, \quad
\sum \|x_k\|_2\x \|y_k\|_2<\infty \right\}.
\end{equation}
By (\ref{invo}), 
$$
X\subset \bar\theta z H^1\cap\theta\overline{zH^1}
$$
and
$$
X_a=\{\bar z\theta f: f\in X\}\subset 
H^1\cap\bar z\theta^2H^1_-\subset K_{\theta^2}^1.
$$
The norms in these spaces are defined as infimum of 
$\sum \|x_k\|_2\x \|y_k\|_2$ over all possible representations, thus
$X$, $X_a$ are Banach spaces.

\begin{Prop}
\label{classx}
$1)$ Any nonnegative element of $X$ can be written as $|g|^2$, $g\in\Kth$.

$2)$ Any element of $X$ can be represented as a linear combination of four
nonnegative elements of $X$.

$3)$ Every element of $X, X_a$ admits a representation as a sum containing
only four summands in the definition of these spaces, and the norm of
each summand in the space does not exceed the norm of the initial element 
of the space.

\end{Prop}

\beginpf $1)$ Let $f =\sum x_k\bar y_k \in X$, $f\ge 0$. 
Since $\bar z\theta\bar y_k \in \Kth$, we have $\bar z\theta f \in H^1$.
Then, by Dyakonov's result \cite{dya}, $f=|g|^2$ for some $g \in \Kth$
(proof: take the outer function with modulus $f^{1/2}$ on $\T$ as $g$; 
then $\bar z\theta\bar g \in H^2$ and hence $g\in \Kth$). 

$2)$ Since $X$ is symmetric with respect to complex conjugation,
it suffices to show that real functions from $X$ may be represented
as a difference of two nonnegative functions from $X$. 
The real part of a function from $X$ of the form  $\sum x_k\bar y_k$
with $x_k, y_k\in\Kth$, $\sum \|x_k\|_2\x\|y_k\|_2<\infty$, is
$$
\frac{1}{2}\sum (x_k\bar y_k+\bar x_k y_k)
=\sum\left|\frac{x_k+y_k}{2}\right|^2-\sum\left|\frac{x_k-y_k}{2}\right|^2,
$$
which is the desired representation. We may suppose that $\|x_k\|=\|y_k\|$ 
for every $k$, then each of the norms 
$\left\Vert \sum\left|\frac{x_k+y_k}{2}\right|^2 \right\Vert_X$,
$\left\Vert \sum\left|\frac{x_k-y_k}{2}\right|^2 \right\Vert_X$
obviously does not exceed $\sum \|x_k\|_2\x\|y_k\|_2$.

$3)$ For the space $X$ this directly follows from $1)$ and $2)$, 
for $X_a$ it remains to use the relation $X_a=\bar z\theta X$,
which is a consequence of (\ref{invo}). \qed

\bigskip

Given a function $f$ in the unit disk, define functions $f_r$, $0<r<1$, by 
$f_r(z)=f(rz)$. We may think of functions $f\in X_a$ as analytic functions in 
$\D$. For $f\in X_a$ write $f=\sum x_k y_k$ with $x_k, y_k\in\Kth$ and
$\sum \|x_k\|_2\x \|y_k\|_2 < \infty$. We have $f_r=\sum (x_k)_r (y_k)_r$.
Now  it follows from Proposition~\ref{polt} that, for any 
$\mu\in \C_2^+(\theta)$, the embedding of the space $X_a$ into $L^1(\mu)$ 
is a well-defined bounded map realized by the limit of $f_r$ in $L^1(\mu)$
as $r\nearrow 1$. 

\medskip

We will need the following important lemma.

\begin{Lem} 
\label{uff}
Let $\mu\in\C_2(\theta)$ and let $x_k, y_k\in\Kth$, 
$\sum \|x_k\|_2\x \|y_k\|_2<\infty$. 
If $\sum x_k\bar y_k=0$ in the space $X$, then also $\sum x_k\bar y_k=0$ \, 
$|\mu|$-almost everywhere.
\end{Lem}

In other words, the embedding $X\into L^1(|\mu|)$ is well defined. 

\beginpf
There is no loss of generality if we assume that $\mu\in\C_2^+(\theta)$.
As above, let $J$ stand for the embedding $\Kth\into L^2(\mu)$.
If $g\in\Kth$, then $\tilde g\in\Kth$, where $\tilde g=\bar z\theta\bar g$. 
By Proposition~\ref{polt} the functions $\tilde g_r$ have a limit as
$r\nearrow 1$, and we want to show that
\begin{equation}
\label{invol}
\lim_{r\to 1-}\tilde g_r = \bar z\theta\bar g \qquad {\rm in} \; L^2(\mu).
\end{equation} 
It suffices to check this relation on a dense set. It is easily seen
that for reproducing kernels 
$\frac{1-\overline{\theta(\lambda)\theta}}{1-\bar\lambda z}$
this property is equivalent to the fact that $|\theta_r|^2\to 1$
proved in Proposition~\ref{mod1}.

Take $x_k, y_k\in\Kth$ such that $\sum \|x_k\|_2\x \|y_k\|_2<\infty$
and $\sum x_k\bar y_k=0$. Consider the functions $\tilde y_k\in\Kth$,
$\tilde y_k=\bar z\theta\bar y_k$. By (\ref{invol}) we have 
$(\tilde y_k)_r\to \bar z\theta\bar y_k$ in $L^2(\mu)$.
The formula $\sum x_k\tilde y_k$ determines the zero element of $X_a$, 
hence $\sum (x_k)_r(\tilde y_k)_r=(\sum x_k\tilde y_k)_r=0$.
We obtain 
$$
\sum x_k\cdot\bar z\theta\bar y_k= \lim_{r\to 1-} \sum (x_k)_r(\tilde y_k)_r=0
$$
in the norm of the space $L^1(\mu)$.
Since $\theta\ne 0$ $\mu$-a.e. (e.g., by Proposition~\ref{mod1}), 
we conclude that $\sum x_k\bar y_k=0$\, $\mu$-almost everywhere. \qed
\bigskip
\medskip

\section{Proofs of Theorems \ref{t1}--\ref{t2}}\label{proofs}

\noindent {\bf Proof of Theorem \ref{t1}.} \;
1) Let $A$ be a nonnegative bounded truncated Toeplitz operator with symbol
$\phi$. Denote by $X_c$ the set of all continuous functions from $X$ and
define the functional $l$ on $X_c$ by $l: f\mapsto\int\phi f\, dm$.
If $f\in X_c$, $f\ge 0$, then $lf\ge 0$. Indeed, by item 1) of
Proposition~\ref{classx}, there exists a function $g\in\Kth$ such that 
$|g|^2=f$ ($g$ turns out to be bounded), and hence 
$$
lf=\int\phi f\,dm=\int\phi g\bar g\, dm=(Ag, g)\ge 0.
$$

Assume first that $\theta(0)=0$, then $1\in X$. Extend the functional $l$ 
to the space $C(\T)$ of all continuous functions on $\T$ by 
the Hahn--Banach theorem so that the norm of the extended functional
equal the norm of $l$. Since $1\in X$, it will be nonnegative automatically, 
hence $lf=\int f\, d\mu$, $f\in X_c$, for some nonnegative Borel measure $\mu$ 
on $\T$. The map $\Kth\to L^2(\mu)$, which takes 
continuous functions to their traces on the support of $\mu$, 
is bounded. Indeed, if $g\in\Kth$ is continuous, then
$$
\int |g|^2d\mu=(Ag, g)\le \|A\|\cdot \|g\|_2^2.
$$
This proves that $\mu\in\C_2^+(\theta)$. 
By linearity and continuity the relation $\int |g|^2d\mu=(Ag, g)$, $g\in\Kth$,
implies $\int x\bar y\, d\mu=(Ax, y)$ for all $x, y\in\Kth$, hence $A=A_\mu$.

If $w = \theta(0)\ne 0$, consider so-called Crofoot's transform 
\begin{equation}
\label{crof}
U: f \mapsto \sqrt{1 - |w|^2} \frac{f}{1 - \overline{w}\theta}
\end{equation}
which is a unitary map of  $K_\theta$ onto $K_{\Theta}$, 
where $\Theta = \frac{\theta - w}{1 - \overline{w}\theta}$
is the Frostman shift of $\theta$. 
Take a bounded truncated Toeplitz operator 
$A \geq 0$ acting on the space $K_\theta$. 
By \cite[Theorem 13.2]{sar} the operator 
$UAU^* \geq 0$ is a bounded truncated Toeplitz 
operator on $K_\Theta$. Note that $\Theta(0) = 0$. 
Let $\mu \in \C^+_2(\Theta)$ be a quasisymbol of $UAU^*$. 
Then the measure $\nu = \frac{1 - |w|^2}{|1 - \overline{w}\theta|^2} \mu$ 
is a quasisymbol of $A$. 
Indeed, $\nu\in\C^+_2(\theta)$, and from (\ref{crof}) it follows that 
for any $f,g \in K_\theta$, we have
$$ 
(Af,g)=(UAU^*Uf,Ug)=\int Uf\cdot\overline{Ug}\,d\mu=\int f\bar{g} \,d\nu.
$$
Thus, each nonnegative bounded truncated Toeplitz operator admits 
a quasisymbol from $\C_2^+(\theta)$.    

2) Let $A$ be a bounded truncated Toeplitz operator. It may be represented 
in the form $A=A_1-A_2+iA_3-iA_4$, where all $A_i$, $i=1, 2, 3, 4$, are 
nonnegative truncated Toeplitz operators. Indeed, $A^\ast$ is a truncated
Toeplitz operators as well, which allows us to consider only selfadjoint
operators. The identity operator $I$ is trivially a truncated Toeplitz
operator (with symbol 1), and $A$ is the difference of two nonnegative
operators $\|A\|\cdot I$ and $\|A\|\cdot I-A$. For each $A_i$ construct
$\mu_i$ as above. It remains to take $\mu=\mu_1-\mu_2+i\mu_3-i\mu_4$. \qed

\medskip

\noindent {\bf Proof of Theorem \ref{c1}.} \;
If $A$ has a bounded symbol $\phi$, then $A=A_\mu$ with $d\mu = \phi dm$, 
and $\mu\in\C_1(\theta^2)$. 

Now let $\mu\in\C_1(\theta^2)$. We need to prove that $A_\mu$ coincides with 
a truncated Toeplitz operator with a bounded symbol. Define the functional
$l: f\mapsto \int f\, d\mu$  on functions from $X$ which are finite 
sums of the functions of the form $x_k\bar y_k$ with 
$x_k, y_k\in\Kth$. For $f=\sum x_k\bar y_k$ we have 
$$
\left|\int f\, d\mu\right|\le \int |f| d|\mu| \le C \|f\|_1,
$$
since $\theta \bar z f \in K_{\theta^2}^1$ 
and $\mu\in\C_1(\theta^2)$.
Hence, the functional $l$ can be continuously extended to $L^1$, and so
there exists a function $\phi\in L^\infty$ such that 
$l(f) = \int \phi f\, dm$, $f\in X$.
Hence, for any $x,y \in K_\theta^2$, we have
$$
\int x\bar y\,d\mu =l(x\bar y) = \int\phi\,x\bar y\, dm= (A_\phi x, y),
$$
and thus $A_\mu = A_\phi$.  \qed

\medskip

\noindent {\bf Proof of Theorem \ref{t2}.} \;
$1)$ First, we verify that the functional (\ref{overx}) is well defined
for any operator $A\in\Tau(\theta)$. We need to prove that $\sum (Ax_k, y_k)=0$ 
if $x_k, y_k\in K_\theta$, $\sum \|x_k\|_2\x \|y_k\|_2<\infty$, 
and $\sum x_k \bar y_k=0$ almost everywhere with respect to the Lebesgue 
measure. To this end, apply Theorem \ref{t1} and find a measure 
$\mu\in\C_2(\theta)$ such that $(Ax, y)=\int x \bar y\,d\mu$ for 
all $x,y \in K_\theta$. Lemma \ref{uff} 
holds for all complex measures from $\C_2(\theta)$; we conclude that
$\sum x_k \bar y_k=0$ \, $\mu$-almost everywhere. 
By the definition of $A_\mu$ we have
$$
\sum (Ax_k, y_k) = \int\left(\sum x_k \bar y_k\right)d\mu.
$$
Thus $\sum (Ax_k, y_k)=0$, and the functional is defined correctly. 
  
Now prove the equality $\|\Phi_A\| = \|A\|.$ Indeed, for any function 
$\sum x_k \bar{y}_k \in X$ we have 
$$
\left|\Phi_A\left(\sum x_k \bar{y}_k\right)\right| = 
\left|\sum (Ax_k, y_k)\right| \le \|A \|  \sum \|x_k\|_2\x \|y_k\|_2.
$$
Hence $\|\Phi_A\| \le \|A\|$. On the other hand, 
for any unit norm vectors $x,y \in K^{2}_{\theta}$ 
we have $\|x\bar{y}\|_X \le 1$  and 
$$
\| A \| =\sup_{\|x\|_2, \|y\|_2\le 1} |(Ax,y)| = \sup_{\|x\|_2, \|y\|_2\le 1}
|\Phi_A(x\bar{y})| \le \|\Phi_A\|.
$$ 
This proves the inverse inequality.

It remains to show that any linear continuous functional $\Phi$ on $X$ 
may be represented in the form $\Phi = \Phi_A$ for some (unique) 
truncated Toeplitz operator $A$. Take a continuous functional $\Phi$ on $X$ 
and define the operator $A_\Phi$ by its bilinear form:
$(A_\Phi x, y)\stackrel{\mathrm{def}}{=} \Phi(x\bar{y})$.
If $f, zf\in\Kth$, we have
$$
(A_\Phi f,f) = \Phi(|f|^2) = \Phi(|zf|^2) = (A_\Phi zf,zf)
$$ 
Now, applying Theorem \ref{sarason}, we obtain $A\in\Tau(\theta)$. The uniqueness 
of $A$ is a consequence of the relation $\|A_\Phi\| = \|\Phi\|$. 

$2)$ Consider the duality $\langle f, A\rangle = \sum (Ax_k, y_k)$
where $A\in \Tau(\theta)$, $f=\sum x_k\bar y_k \in X$ 
(by Lemma \ref{uff} $\langle f, A\rangle$ does not depend on the 
choice of factorization).
We need to prove that every continuous functional $\Phi$ 
on the space $\Tau_0(\theta)$ of 
all compact truncated Toeplitz operators is realized by an element of $X$. 
Extend $\Phi$ by Hahn--Banach theorem to the space
of all compact operators in $\Kth$. The trace class is the dual space
to the class of all compact operators, hence the functional may be 
written in the form $\Phi(A) = \sum (Ax_k, y_k)$, 
$A \in \Tau_0(\theta)$, 
for some $x_k, y_k\in\Kth$ with $\|\Phi\| = \sum \|x_k\|_2\x\|y_k\|_2$.
Then $f=\sum x_k\bar y_k \in X$ and
$\Phi(A) = \langle f, A\rangle$. Repeating the arguments from 1)
we conclude that $\|\Phi\| \le \|f\|_X$. On the other hand, 
$\|f\|_X \le \sum \|x_k\|_2\x\|y_k\|_2 = \|\Phi\|$.  \qed

\medskip

Now we prove two corollaries which give us additional 
information on the structure of the space 
of truncated Toeplitz operators.

\smallskip

\begin{Cor}\label{comp}
The closed (in the norm) linear span of rank-one 
truncated Toeplitz operators
coincides with the set of all compact truncated Toeplitz operators.
\end{Cor}

\beginpf
For $\lambda\in \D$, denote by $k_\lambda, \tilde k_\lambda$ 
the functions from $\Kth$, 
$$
k_\lambda(z)=\frac{1-\overline{\theta(\lambda)}\theta(z)}{1-\bar\lambda z},
\qquad
\tilde k_\lambda(z)=\frac{\theta(z)-\theta(\lambda)}{z-\lambda}
$$
(recall that $k_\lambda$ is the reproducing kernel for the space $\Kth$). 
If $x, y\in\Kth$, then $(x, k_\lambda)=x(\lambda)$, 
$(\tilde y, \tilde k_\lambda)=\overline{y(\lambda)}$. 
It is shown in \cite{sar} that the operators 
$\Alambda=(\cdot, k_\lambda)\tilde k_\lambda$
are rank-one truncated Toeplitz operators. Take $f\in X$ as an element of
the dual space to the class of all compact truncated Toeplitz operators.
Set $g=\bar z\theta f\in X_a$, let $g=\sum x_k y_k$ with $x_k, y_k\in\Kth$.
The following formula illustrates the duality on rank-one truncated Toeplitz 
operators:
$$
\begin{aligned}
\langle f, \Alambda \rangle=&\big\langle
\sum x_k\cdot\overline{\tilde y_k}, \Alambda\big\rangle
=\sum (\Alambda x_k, \tilde y_k)\\
=&\sum x_k(\lambda)\cdot (\tilde k_\lambda, \tilde y_k)
=\sum x_k(\lambda)\cdot y_k(\lambda)=g(\lambda).
\end{aligned}
$$
Suppose that $f$ annihilates all operators of the form $\Alambda$
with $|\lambda|<1$. Then $g\equiv 0$ in $\D$, hence $f$ is
the zero element of $X$. 
\qed

\medskip

The space of all bounded linear operators on a Hilbert space is dual to
the space of trace class operators. This duality generates the
ultraweak topology on the former space. Formally, the ultraweak topology is
stronger than the weak operator topology, but on the subspace of all
truncated Toeplitz operators they coincide. 

\begin{Cor}\label{WOT}
The weak operator topology on $\Tau$ coincides with the ultraweak topology.  
\end{Cor}

\beginpf Any ultraweakly continuous functional $\Phi$ on $\Tau$ is
generated by some trace class operator $\sum_k (\cdot, y_k) x_k$,
where $x_k,y_k \in \Kth$, $\sum_k \|x_k\|_2 \x \|y_k\|_2 < \infty$,
and is of the form
$$
\Phi(A) = \sum_k (A x_k, y_k).
$$ 
The function $f = \sum_k x_k \bar y_k $ belongs to the space $X$. 
It follows from Proposition \ref{classx} that there exist 
$f_1, g_1 \ldots f_4,g_4 \in \Kth$ such that 
$f = \sum_{k=1}^{4} f_k \bar g_k$. Therefore, 
by the duality from Theorem \ref{t2},
$$
\Phi(A) = \langle f, A \rangle = 
\big\langle \sum_{k=1}^4 f_k \bar g_k, A \big\rangle 
= (A f_1,g_1) + \ldots + (A f_4, g_4).
$$ 
Now the statement of the corollary is obvious.  
\qed

\bigskip


\section{Proof of Theorem \ref{t3}}

Throughout this section we will assume, for simplicity, 
that $\theta(0) = 0$. The general case follows immediately 
by means of the transform (\ref{crof}) (note that in this case
$\C_p(\Theta) = \C_p(\theta)$ for any $p$, see, e.g., \cite[Theorem 1.1]{al5}).
\smallskip

For the proof of Theorem \ref{t3} we need the following obvious lemma 
(see \cite{al5}) based on the relations $K_{\theta^2}=\Kth\oplus\theta\Kth$
and $\Kth\cdot\Kth \subset K_{\theta^2}^1$.

\begin{Lem} 
\label{lemC(th)=C(th^2)}
For any inner function $\theta$ we have
$\C_2(\theta)=\C_2(\theta^2)$ and $\C_1(\theta^2)\subset\C_2(\theta^2)$.
If $\theta(0)=0$, we also have $\C_p(\theta^2)=\C_p(\theta^2/z)$
for any $p$. The same equalities or inclusions hold for the classes
$\Dd_p(\theta)$.
\end{Lem}

\noindent {\bf Proof of Theorem \ref{t3}.} \; 
$3) \Rightarrow 2)$. We will establish
condition $2')$, which is formally stronger than $2)$. By Lemma \ref{lemC(th)=C(th^2)} it suffices to prove 
the inclusion $\Dd_2(\theta) \subset \Dd_1(\theta^2/z)$. Take a complex 
measure $\mu\in\Dd_2(\theta)$. We must check that the embedding 
$K^{1}_{\theta^2/z} \into L^1(|\mu|)$ is a bounded operator. 
By condition~$3)$, 
there is a positive constant $c_1$ such that any function 
$f \in K^{1}_{\theta^2/z} $ can be represented in the form 
$f = \sum_{k=1}^{\infty} f_k g_k$, where the functions 
$f_k, g_k$ are in $\Kth$ and
$\sum_{k=1}^{\infty} \|f_k\|_2\x\|g_k\|_2 
\le c_1 \|f\|_1$.
Since $\mu \in \Dd_2(\theta)$, we have 
$\sum_{k = 1}^{\infty} \|f_k\|_{L^2(|\mu|)}\x \|g_k\|_{L^2(|\mu|)}<\infty$,
so the series converges in $L^1(|\mu|)$ and $f\in L^1(|\mu|)$.
Moreover, 
$$
\begin{aligned}
\int |f|\,d|\mu| 
& =
\int\left|\sum_{k=1}^\infty f_k g_k\right|\,d|\mu| \le 
\sum_{k = 1}^{\infty} \|f_k\|_{L^2(|\mu|)}\x \|g_k\|_{L^2(|\mu|)} \\
&\le 
c_2 \sum_{k = 1}^{\infty} \|f_k\|_2\x \|g_k\|_2 \le
c_1 c_2 \|f\|_1,
\end{aligned}
$$
and, thus, $\mu \in \Dd_1(\theta^2/z)$. 
Therefore, we have $\Dd_2(\theta)\subset \Dd_1(\theta^2/z)$
which implies $\Dd_2(\theta^2)=\Dd_1(\theta^2)$
(and, in particular, $\C_2(\theta^2)=\C_1(\theta^2)$).    

\medskip

The implication $2)\Rightarrow 1)$ follows directly from 
Theorems \ref{t1} and \ref{c1} (even with the weaker
assumption $\C_1(\theta^2)=\C_2(\theta^2)$).

\medskip

$1) \Rightarrow 3)$. Condition $3)$ can 
be written in the form $X_a = K^{1}_{\theta^2/z}$ or, equivalently, 
as $X = z\bar{\theta} K^{1}_{\theta^2/z}$ (see Section 4). 
In the general case we have that $X$ is a dense subset of 
$z\bar{\theta} K^{1}_{\theta^2/z}$. By the Closed Graph theorem,
$X = z\bar{\theta} K^{1}_{\theta^2/z}$ if and only if the norms in 
the spaces $X$ and $K^{1}_{\theta^2/z}$ are equivalent. Take an arbitrary 
function $h = \sum x_k \bar{y}_k \in X$. Clearly, we have 
$\|h\|_1 \le \|h \|_X$.  On the other hand, it follows from Theorem \ref{t2}
that
\begin{equation}
\label{dual}
\|h\|_X = \sup \left\{\left|\sum (Ax_k,y_k)\right|: \, A \in \Tau(\theta), \|A\| 
\le 1  \right\}.
\end{equation}   
The Closed Graph theorem and condition $1)$ guarantee the existence of 
a bounded symbol $f_A \in L^\infty$ for any operator $A \in \Tau(\theta)$ with 
$\|f_A \|_\infty \le c \|A\|$. Therefore, the supremum in (\ref{dual})
does not exceed
$$
\sup\left\{\left|\sum (fx_k,y_k)\right|: f \in L^\infty, \|f\|_\infty 
\le c \right\} = c\! \sup_{\|f\|_\infty \le 1}\left|\int\! f\x
\sum x_k \bar{y}_k \,dm\right|.
$$ 
Thus, $\|h\|_X \le c \|h\|_1$, $h\in X$, which proves the theorem.  \qed

\bigskip
\section{One-component inner functions}\label{onecomp}

As we have noted in Section \ref{main}, one-component inner functions
satisfy the condition 2) of Theorem \ref{t3}: $\C_2(\theta)=\C_1(\theta^2)$.
It is also possible to show that one-component inner functions satisfy
the factorization condition 1) in Theorem \ref{t3}.
We will start with a particular case of the Paley--Wiener spaces.
 
\begin{Ex} 
{\rm Let $\Theta_a(z) = \exp(iaz) $, $a>0$, 
be an inner function in the upper half-plane.
Then for the corresponding model subspace 
we have $ K_{\Theta_a}^p = \PW_a^p \cap H^p$.
Note, that the model subspaces in the half-plane case are defined as
$K_\Theta^p = H^p \cap \Theta \overline{H^p}$, the involution is given by
$f\mapsto \Theta\bar f$, and so $fg\in K_{\Theta^2}^1$ for any $f,g \in K^2_\Theta$.
Thus, in view of Proposition \ref{classx}, the factorization for 
the corresponding space $X$ is equivalent to the following 
property:  for any $f\in \PW^1_{2a}$ which takes real values there exist $g\in 
\PW^1_{2a}$ such that $|f|\le g$. This can be easily achieved.  Let $a=\pi/2$.
Put
$$
g(z)= \sum_{n\in\Z} c_n \frac{\sin^2 \frac{\pi}{2}(t-n)}{(t-n)^2},
$$ 
where $c_n = \max_{[n,n+1]} |f|$. By the 
Plancherel--P\'olya inequality (see, e.g., \cite[Lecture 20]{levin}), 
$\sum_n c_n \le C\|f\|_1$, and so 
$g\in \PW^1_{\pi}$. Also, if $t\in [n, n+1]$, then  
$$
|f(t)| \le c_n \le  c_n \frac{\sin^2 \frac{\pi}{2}(t-n)}{(t-n)^2} \le g(t).
$$   }
\qed
\end{Ex}

An analogous argument works for general one-component 
inner function. Let $\theta$ be an inner function in the unit disk. 
In view of Proposition \ref{classx}, 
the property 1) in Theorem 
\ref{t3} will be obtained as soon as we prove the following theorem:

\begin{Thm} 
\label{onecomp2}
For any real-valued element $f$ of $X$ there exist $t_n\in \tz$ and 
$c_n>0$ such that 
$$
g = \sum c_n |k_{t_n}|^2 \in X,
$$ 
$\|g\|_1 \le C\|f\|_1$, and $|f| \le g$ on $\mathbb{R}$.
\end{Thm}

First we collect some known properties of one-component inner functions.
\smallskip

(i) Let $\rho(\theta)$ be the so-called {\it spectrum} of the inner
function $\theta$, that is, the set of all 
$\zeta\in\overline {\mathbb{D}}$
such that $\liminf\limits_{z\to\zeta,\, z\in\mathbb{D}}|\theta(z)|=0$.
Then $\theta$, as well as any element of $\kp$, 
has an analytic extension across any subarc  
of the set $\mathbb{T}\setminus \sigma(\theta)$.

It is shown in \cite{al2} that for a one-component inner function 
$\sigma_\alpha(\rho(\theta))=0 $ for any Clark measure $\sigma_\alpha$
defined by (\ref{clark}).
Thus, all Clark measures are purely atomic and supported on the set
$\T\setminus \rho(\theta)$.
\medskip

(ii) On each arc of the set $\T\setminus \rho(\theta)$,
there exists a smooth increasing branch of the argument of $\theta$
(denote it by $\psi$) and the change of the argument between
two neighboring points $t_n$ and $t_{n+1}$
from the support of one Clark measure
is exactly $2\pi$.
\medskip

(iii) By $(t_n, t_{n+1})$ we denote the closed arc with endpoints 
$t_n, t_{n+1}$, which contains no other points from the Clark measure support.
There exists a constant $A = A(\theta)$ such that for any two points
$t_n$ and $t_{n+1}$ satisfying $|\psi(t_{n+1}) -\psi(t_n)|=2 \pi$
and for any $s, t$ from the arc $(t_n, t_{n+1})$,
\begin{equation}
\label{3}
A^{-1}  \le \frac{|\theta'(s)|}{|\theta'(t)|} \le A,
\end{equation}
that is, $|\theta'|$ is almost constant, when the 
change of the argument is small.
This follows from the results of \cite{al2}, 
a detailed proof may be found in \cite[Lemma 5.1]{bardya}.
\medskip

(iv) If $\theta$ is one-component, then 
${\mathcal C}_1(\theta) = {\mathcal C}_2 (\theta)$. 
The same holds for the function $\theta^2$ which is also one-component.
By Lemma \ref{lemC(th)=C(th^2)}, 
${\mathcal C}_1(\theta^2)= {\mathcal C}_2(\theta^2) = {\mathcal C}_2 (\theta)$,
and so there is a constant $B$ such that for any measure in 
${\mathcal C}^+_2(\theta)$,
\begin{equation}
\label{3a}
\sup_{f\in K^1_{\theta^2}} \frac{\|f\|_{L^1(\mu)}}{\|f\|_1} \le
B \sup_{f\in K^2_{\theta}} \frac{\|f\|_{L^2(\mu)}}{\|f\|_2}.
\end{equation}

(v) Let $\{t_n\}$  be the support of some Clark measure for $\theta$
and let $s_n \in (t_n, t_{n+1})$. There exists a constant $C=C(\theta)$
which does not depend on $\{t_n\}$ and $\{s_n\}$ such that
for any $f\in K_{\theta}^2$,
$$
\sum_n \frac{|f(s_n)|^2}{|\theta'(s_n)|} \le C\|f\|_2^2.
$$
This follows from the stability result due to Cohn \cite[Theorem 3]{cohn86}.
\medskip 

So (v) means that for the measures of the form 
$\sum_n |\theta'(s_n)|^{-1} \delta_{s_n}$ the supremum in the 
right-hand side of (\ref{3a}) is uniformly bounded.
From this, (iii) and (iv) we have the following Plancherel--Polya type inequality:

\begin{Cor}
\label{carl1}
Let $\{t_n\}$ be the support of some Clark measure for $\theta$
and let $s_n, u_n \in (t_n, t_{n+1})$. There exists a constant $C=C(\theta)$
which does not depend on $\{t_n\}$, $\{s_n\}$, $\{u_n\}$, such that
for any $f\in X\subset K_{\theta^2}^1$,
\begin{equation}
\label{4}
\sum_n \frac{|f(s_n)|}{|\theta'(u_n)|} \le C\|f\|_1.
\end{equation}
\end{Cor}

{\bf Proof of Theorem~\ref{onecomp2}.}
Let $f\in X$. Take two Clark bases corresponding to $1$ and $-1$,
and let $\{t_n\}$ be the union of their supports. 
If $t_n$ and $t_{n+1}$ are two neighbor points from our set, then
$$
\int_{(t_n, t_{n+1})} |\theta'(t)| dm(t) = \pi.
$$
If we write $t_n = e^{ix_n}$ and take the branch of the argument
$\psi$ so that $\theta (e^{ix}) = e^{2i\psi(x)}$,
then $|\psi(x_{n+1})-\psi(x_n)| = \pi/2$.

Let $c_n = \sup_{t\in (t_n, t_{n+1})} |f(t)|$ 
and put, for some constant $D$ whose value will be specified later, 
$$
g(z)  = D\sum  c_n \frac{|k_{t_n}(z)|^2}{|\theta'(t_n)|^2}.
$$
Then $g \in L^1$ since the series converges in $L^1$-norm. Indeed, 
$c_n = |f(s_n)|$ for some $s_n \in (t_n, t_{n+1})$ and 
$$
\sum |c_n| \frac{\|k_{t_n}^2\|_1}{|\theta'(t_n)|^2} = 
\sum \frac{|f(s_n)|}{|\theta'(t_n)|} \le C\|f\|_1
$$
by Corollary \ref{carl1}. Also $g\in X$ and $g\ge 0$. 

It remains to show that $g\ge |f|$. Let $t = e^{ix}\in (t_n, t_{n+1})$.
We have
\begin{equation}
\label{5}
|k_{t_n}(t)|= \bigg| \frac{\theta(t) -\theta(t_n)}{t-t_n}\bigg| 
= \bigg|2\frac{\sin (\psi(x)-\psi(x_n))}
{e^{ix} -e^{ix_n}}\bigg|.
\end{equation}
Since $|\psi(x)-\psi(x_n)| \le \pi/2$,
we have $|\sin (\psi(x)-\psi(x_n))| \ge 2 |\psi(x)-\psi(x_n)|/\pi$.
Of course we have $|e^{ix} -e^{ix_n}|\le |x-x_n|$. Hence,
the last quantity in (\ref{5}) is
$$
\ge \frac{4}{\pi}\cdot
\bigg|\frac{\psi(x)-\psi(x_n)}{x -x_n} \bigg| = \frac{4\psi'(y_n)}{\pi}
$$
for some $y_n\in [x_n, x]$. If we put $u_n = e^{iy_n}$
we have $\psi'(y_n) = |\theta'(u_n)|/2$.
Thus, we have shown that $|k_{t_n}(t)| \ge 2|\theta'(u_n)| /\pi$
for some $u_n \in (t_n, t_{n+1})$.
Hence, if we take $D> \pi^2 A^2/4$, then
$$
g(t) > Dc_n \frac{|k_{t_n}(t)|^2}{|\theta'(t_n)|^2}
\ge D\frac{4|\theta'(u_n)|^2}{\pi^2 |\theta'(t_n)|^2} c_n
\ge \frac{4D}{\pi^2 A^2} c_n \ge c_n \ge |f(t)|.  
$$
The theorem is proved. \qed

\enddocument